\documentclass[12pt]{article}
\usepackage[T2A]{fontenc}
\usepackage[cp1251]{inputenc}
\usepackage[english]{babel}
\usepackage{amsmath,latexsym,amsthm,amsfonts,tipa,upgreek}
\usepackage{amssymb,amscd,graphpap}

\newcommand\NN{{\mathbb N}}
\newcommand\ZZ{{\mathbb Z}}

\newcommand\w{{\omega}}

\newcommand\PP{{\mathcal P}}

\newcommand\vt{{\Updelta}}
\newtheorem{Th}{Theorem}

\newtheorem{Ps}{Proposition}

\begin{document}

\title{Scattered Subsets of Groups}
\author{T.O.~Banakh (Lviv Nat. Ivan Franko Univ.)\\I.V.~Protasov, S.V.~Slobodianiuk (Kyiv Nat. Taras Shevchenko Univ.)}
\date{}
\maketitle

\begin{abstract}
We define the scattered subsets of a group as asymptotic counterparts of scattered subspaces of a topological space, and prove that a subset $A$ of a group $G$ is scattered if and only if $A$ contains no piecewise shifted $IP$-subsets. For an amenable group $G$ and a scattered subspace $A$ of $G$, we show that $\mu(A)=0$ for each left invariant Banach measure $\mu$ on $G$.

\end{abstract}

\section{Introduction}
Given a discrete space $X$, we take the points of $\beta X$, the Stone-$\check{C}$ech compactificationof $X$, to be the ultrafilters on $X$, with the points of $X$ identified with the principal ultrafilters on $X$. The topology on $\beta X$ can be defined by stating that the sets of the form $\overline{A}=\{p\in\beta X: A\in p\}$, where $A$ is a subset of $X$, form a base for the open sets. We note that the sets of this form are clopen and that, for any $p\in\beta X$ and $A\subseteq X$, $A\in p$ if and only if $p\in\overline{A}$. For any $A\subseteq X$ we denote $A^*=\overline{A}\cap G^*$, where $G^*=\beta G\setminus G$. The universal property of $\beta G$ states that every mapping $f:X\to Y$, where $Y$ is a compact Hausdorff space, can be extended to the continuous mapping $f^\beta:\beta X\to X$.

Now let $G$ be a discrete group. Using the universal property of $\beta G$, we can extend the group multiplication from $G$ to $\beta G$ in two steps. Given $g\in G$, the mapping $$x\mapsto gx:\text{ }G\to\beta G$$
extends to the continuous mapping $$q\mapsto gq:\text{ }\beta G\to\beta G.$$
Then, for each $q\in\beta G$, we extend the mapping  $g\mapsto gq$ defined from $G$ into $\beta G$ to the continuous mapping $$p\mapsto pq:\text{ }\beta G\to\beta G.$$
Te product $pq$ of the ultrafilters $p$, $q$ can also be defined by the rule: given a subset $A\subseteq G$,
$$A\in pq\leftrightarrow \{g\in G:g^{-1}A\in q\}\in p.$$
To describe the base for $pq$, we take any element $P\in p$ and, for every $x\in P$, choose some element $Q_x\in q$. Then $\cup_{x\in p}xQ_x\in pq$, and the family of subsets of this form is a base for the ultrafilter $pq$.

By the construction, the binary operation $(p,q)\mapsto pq$ is associative, so $\beta G$ is a semigroup, and $G^*$ is a subsemigroup of $\beta G$. For each $q\in\beta G$, the right shift $x\mapsto xq$ is continuous, and the left shift $x\mapsto xq$ is continuous for each $g\in G$.

For the structure of a compact right topological semigroup $\beta G$ and plenty of its applications to combinatorics, topological algebra and functional analysis see \cite{b1}, \cite{b2}, \cite{b3}, \cite{b4}, \cite{b5}.

Given a subset $A$ of a group $G$ and an ultrafilter $p\in G^*$, we define a {\em $p$-companion} of $A$ by
$$\vt_p(A)=A^*\cap Gp=\{gp: g\in G, A\in gp\},$$
and say that a subset $S$ of $G^*$ is an {\em ultracompanion} of $A$ if $S=\vt_p(A)$ for some $p\in G^*$. For ultracompanions of subsets of groups and metric spaces see \cite{b6}, \cite{b7}.

Clearly, $A$ is finite if and only if $\vt_p(A)=\varnothing$ for each $p\in G^*$.

We say that a subset $A$ of a group $G$ is
\begin{itemize}
\item{} {\em thin} if $|\vt_p(A)|\le1$ for each $p\in G^*$;
\item{} {\em $n$-thin}, $n\in\NN$ if $|\vt_p(A)|\le n$ for each $p\in G^*$;
\item{} {\em sparse} if each ultracompanion of $A$ is finite;
\item{} {\em disparse} if each ultracompanion of $A$ is discrete;
\item{} {\em scattered} if, for each infinite subset $Y$ of $A$, there is $p\in Y^*$ such that $\vt_p(Y)$ is finite.
\end{itemize}

We denote by $[G]^{<\w}$ the family of all finite subsets of $G$. Given any $F\in[G]^{<\w}$ and $g\in G$, we put 
$$B(g,F)=Fg\cup\{g\}$$
and, following \cite{b8}, say that $B(g,F)$ is a {\em ball of radius $F$ around $g$}. For a subset $Y$ of $G$, we put $B_Y(g,F)=Y\cap B(g,F)$. By \cite[Ppoposition 4]{b6}, $Y$ is $n$-thin if and only if for every $F\in[G]^{<\w}$, there exists $H\in[G]^{<\w}$ such that $|B_Y(y,F)|\le n$ for each $y\in Y\setminus H$. For thin subsets of a group, their applications and modifications see \cite{b9}--\cite{b19}.

By \cite[Proposition 5]{b6} and \cite[Theorems 3 and 10]{b20}, for a subset $A$ of a group $G$, the following statements are equivalent
\begin{itemize}
\item[(1)] $A$ is sparse;
\item[(2)] for every infinite subset $X$ of $G$, there exists finite subset $F\subset G$ such that $\bigcap_{g\in F}gA$ is finite;
\item[(3)] for every infinite subset $Y$ of $A$, there exists $F\in[G]^{<\w}$ such that, for every $H\in[G]^{<\w}$, we have $$\{y\in Y: B_A(y,H)\setminus B_A(y,F)=\varnothing\}\neq\varnothing;$$
\item[(4)] $A$ has no subsets asymorphic to the subset $W_2=\{g\in \oplus_\w \ZZ_2: supt g\le2\}$ of the group $\oplus_\w\ZZ_2$, where $supt g$ is the member of non-zero coordinates of $g$.
\end{itemize}
The notion of asymorphisms and coarse equivalence will be defined in the next section. The sparse sets were introdused in \cite{b21} in order to characterise strongly prime ultrafilters in $G^*$, the ultrafilters from $G^*\setminus\overline{G^*G^*}$. More on sparse subsets can be find in \cite{b10}, \cite{b11}, \cite{b16}, \cite{b22}.

In this paper, answering Question 4 from \cite{b6}, we prove that a subset $A$ of a group $G$ is scattered if and only if $A$ is disparse, and characterize the scattered subsets in terms of prohibited subsets. We answer also Question 2 from \cite{b6} proving that each scattered subset of an amenable group is absolute null. The results are exposed in section~\ref{s2}, their proofs in section~\ref{s3}.

\section{Results}~\label{s2}

Our first statement shows that, from the asymptotic point of view \cite{b23}, the scattered subsets of a group can be considered as the counterparts of the scattered subspaces of a topological space.
\begin{Ps}\label{p1} For a subset $A$ of a group $G$, the following two statements are equivalent
\begin{itemize}
\item[(i)] $A$ is scattered;
\item[(ii)] for every infinite subset $Y$ of $A$, there exists $F\in[G]^{<\w}$such that, for every $H\in[G]^{<\w}$, we have $$\{y\in Y: B_Y(y,H)\setminus B_Y(y,F)=\varnothing\}\neq\varnothing.$$
\end{itemize}
\end{Ps}
\begin{Ps}\label{p2}
A subset $A$ of a group $G$ is scattered if and only if, for every countable subgroup $H$ of $G$, $A\cap H$ is scattered in $H$.
\end{Ps}
Let $A$ be a subset of a group $G$, $K\in[G]^{<\w}$. A sequence $a_0,..,a_n$ in $A$ is called {\em $K$-chain} from $a_0$ to $a_n$ if $a_{i+1}\in B(a_i,K)$ for each $i\in\{0,...,n-1\}$. For every $a\in A$, we denote
$$B_A^\square (a,K)=\{b\in A: \text{ there is a }K\text{-chain from } a \text{ to } b\}$$
and, following \cite[Chapter 3]{b24}, say that $A$ is {\em cellular} (or {\em asymptotically zero-dimentional}) if, for every $K\in[G]^{<\w}$, there exists $K'\in[G]^{<\w}$ such that, for each $a\in A$, $$B_A^\square(a,K)\subseteq B_A(a,K').$$
Now we need some more asymptology (see \cite[Chapter 1]{b24}). Let $G$, $H$ be groups, $X\subseteq G$, $Y\subseteq H$. A mapping $f:X\to Y$ is called a {\em $\prec$-mapping} if, for every $F\in[G]^{<\w}$, there exists $K\in[G]^{<\w}$ such that, for every $x\in X$, $$f(B_X(x,F))\subseteq B_Y(f(x),K).$$
If $f$ is a bijection such that $f$ and $f^{-1}$ are $\prec$-mapping, we say that $f$ is an asymorphism. The subsets $X$ and $Y$ are called {\em coarse equivalent} if there exist asymorphic subsets $X'\subseteq X$ and $Y'\subseteq Y$ such that $X\subseteq B_X(X',F)$, $Y\subseteq B_Y(Y',K)$ for some $F\in[G]^{<\w}$ and $K\in[H]^{<\w}$.

Following \cite{b23}, we say, that the set $Y$ of $G$ {\em has no asymptotically isolated balls} if $Y$ does not satisfy Proposition~\ref{p1}$(ii)$: for every $F\in[G]^{<\w}$, there exists $H\in[G]^{<\w}$ such that $B_Y(y,H)\setminus B_Y(y,F)\neq\varnothing$ for each $y\in Y$.

By \cite{b23}, a countable cellular subset $Y$ of $G$ with no asymptotically isolated balls is coarsely equivalent to the group $\oplus_\w\ZZ_2$.
\begin{Ps}\label{p3}
Let $X$ be a countable subset of a group $G$. If $X$ is not cellular then $X$ contains a subset $Y$ coarsely equivalent to $\oplus_\w\ZZ_2$.
\end{Ps}
Let $(g_n)_{n<\w}$ be an injective sequence in a group $G$. The set $$\{g_{i_1}g_{i_2}...g_{i_n}:0\le i_1<i_2<...<i_n<\w\}$$ is called an {\em $IP$-set} \cite[p. 406]{b1}, the abriviation for "infinite dimensional parallelepiped".

Given a sequence $(b_n)_{n<\w}$ in $G$, we say that the set $$\{g_{i_1}g_{i_2}...g_{i_n}b_{i_n}:0\le i_1<i_2<...<i_n<\w\}$$ is a {\em piecewise shifted $IP$-set}.
\begin{Th}\label{t1}
For a subset $A$ of a group $G$, the following statements are equivalent 
\begin{itemize}
\item[(i)] $A$ is scattered;
\item[(ii)] $A$ is disparse;
\item[(iii)] $A$ contains no subsets coarsely equivalent to the group $\oplus_\w\ZZ_2$;
\item[(iv)] $A$ contains no piecewise shifted $IP$-sets.
\end{itemize}
\end{Th}
By the equivalence $(i)\Leftrightarrow(ii)$ and Propositions 10 and 12 from \cite{b6}, the family of all scattered subsets of an infinite group $G$ is a translation invariant ideal in the Boolean algebra of all subsets of $G$ strictly contained in the ideal of all small subsets.

Now we describe some relationships between the left invariant ideals $Sp_G$, $Sc_G$ of all sparse and scattered subsets of a group $G$ on one hand, and closed left ideals of the semigroup $\beta G$.

Let $J$ be a left invariant ideal in the Boolean algebra $\PP_G$ of all subsets of a group $G$. We set $$\hat{J}=\{p\in\beta G: G\setminus A\in p \text{ for each } A\in J\}$$
and note that $\hat{J}$ is a closed left ideal of the semigroup $\beta G$. On the other hand, for a closed left ideal $L$ of $\beta G$, we set 
$$\check{L}=\{A\subseteq G: A\notin p\text{ for each }p\in L\}$$
and note that $\check{L}$ is a left invariant ideal in $\PP_G$. Moreover, $\check{\hat{J}}=J$ and $\hat{\check{L}}=L$.

Clearly, $\hat{[G]^{<\w}}=G^*$ and by Theorem~\ref{t1}, 
$$(\star)\text{ } \hat{Sc_G}=cl\{p\in\beta G: Gp\text{ is discrete in }\beta G\}=$$ 
$$cl\{p\in\beta G : p=\varepsilon p\text{ for some idempotent }\varepsilon\in G^*\}.$$
Given a left invariant ideal $J$ in $\PP_G$ and following \cite{b11}, we define a left invariant ideal $\sigma(J)$ by the rule: $A\in \sigma(J)$ if and only if $\vt_p(A)$ is finite for every $p\in\hat{J}$. Equivalently, $\sigma(J)=\check{cl(G^*\hat{J})}$. Thus, we have $$\hat{Sp_G}=cl(G^*G^*).$$
We say that a left invariant ideal $J$ in $\PP_G$ is {\em sparse-complete} if $\sigma(J)=J$ and denote by $\sigma^*(J)$ the intersection of all sparse-complete ideals containing $J$. Clearly, the sparse-completion $\sigma^*(J)$ is the smallest sparse-complete ideal such that $J\subseteq \sigma^*(J)$. By \cite[Theorem 4(1)]{b11}, $\sigma^*(J)=\bigcup_{n\in\w}\sigma^n(J)$, where $\sigma^0(J)=J$ and $\sigma^{n+1}(J)=\sigma(\sigma^n(J))$. We can prove that $A\in\sigma^n([G]^{<\w})$ if and only if $A$ has no subsets asymorphic to $W_n=\{g\in\oplus_{\w}\ZZ_2:supt g\le n\}$.

By \cite[Theorem 4(2)]{b11}, the ideal $Sp_G$ is not sparse complete. By $(\star)$, the ideal $Sc_G$ is sparse-complete. Hence $\sigma^*([G]^{<\w})\subseteq Sc_G$ but $\sigma^*([G]^{<\w})\neq Sc_G$.

Recall that a subset $A$ of an amenable group $G$ is {\em absolute null} if $\mu(A)=0$ for each left invariant Banach measure $\mu$ on $G$. For sparse subsets, the following theorem was proved in \cite[Theorem 5.1]{b10}.
\begin{Th}\label{t2}
Every scattered subset $A$ of an amenable group $G$ is absolute null.
\end{Th}
Let $A$ be a subset of $\ZZ$. The {\em upper density $\overline{d}(A)$} is denoted by 
$$\overline{d}(A)=\limsup\limits_{n\to\infty}\frac{|A\cap\{-n,-n+1,...,n-1,n\}|}{2n+1}$$
By \cite[Theorem 11.11]{b25}, if $\overline{d}(A)>0$ then $A$ contains a piecewise shifted $IP$-set. We note that Theorem~\ref{t2} generalizes this statement because there exists a Banach measure $\mu$ on $\ZZ$ such that $\overline{d}(A)=\mu(A)$.

In connection with Theorem~\ref{t1}, one may ask if it possible to replace piecewise shifted $IP$-sets to (left or right) shifted $IP$-sets. By Theorem~\ref{t2} and \cite[Theorem 11.6]{b25}, this is impossible.

\section{Proofs}~\label{s3}
{\em Proof of Proposition~\ref{p1}}.

$(i)\Rightarrow (ii)$. We take $p\in Y^*$ such that $\vt_p(Y)$ is finite, so $\vt_p(Y)=Fp$ for some $F\in[G]^{<\w}$. Given any $H\in[G]^{\w}$, we have $hp\notin\vt_p(Y)$ for each $h\in H\setminus F$. Hence $hP_h\cap Y=\varnothing$ for some $P_h\in p$. We put $P=\bigcap_{h\in H\setminus P}P_h$ and note that $$P\subseteq\{y\in Y: B_Y(y,H)\setminus B_Y(y,F)=\varnothing\}.$$
$(ii)\Rightarrow (i)$. We take an infinite subset $Y$ of $A$, choose corresponding $F\in[G]^{<\w}$ and, for each $H\in[G]^{<\w}$, denote $$P_H=\{y\in Y: B_Y(y,H)\setminus B_Y(y,F)=\varnothing\}.$$
By $(ii)$, the family $\{P_H:H\in [G]^{<\w}\}$ has a finite intersection property and $\bigcap_{H\in [G]^{<\w}}P_H=\varnothing$. Hence $\{P_H:H\in [G]^{<\w}\}$ is contained in some ultrafilter $p\in Y^*$. By the choice of $p$, we have $gp\notin\vt_p(Y)$ for each $g\in G\setminus(F\cup\{e\})$, $e$ is the identity of $G$. It follows that $\vt_p(Y)$ is finite so $A$ is scattered.\\\\
{\em Proof of Proposition~\ref{p2}}.

Assume that $A$ is not scattered and choose a subset $Y$ of $A$ which does not satisfy the condition $(ii)$ of Proposition~\ref{p1}. We take an arbitrary $a\in A$ and put $F_0=\{e,a\}$. Then we choose inductively a sequence $(F_n)_{n\in\w}$ in $[G]^{<\w}$ such that
\begin{itemize}
\item[(1)] $F_nF_n^{-1}\subset F_{n+1}$;
\item[(2)] $B_Y(y,F_{n+1})\setminus B_Y(y,F_n)\neq\varnothing$ for every $y\in Y$. 
\end{itemize}
After $\w$ steps, we put $H=\bigcup_{n\in\w}F_n$. By the choice of $F_0$, $Y\cap H\neq\varnothing$. By $(1)$, $H$ is a subgroup. By $(2)$, $(Y\cap H)$ is not scattered in $H$.\\\\
{\em Proof of Proposition~\ref{p3}}.

Replacing $G$ by by the subgroup generating by $X$, we assume that $G$ is countable. We write $G$ as an union of an increasing chain $F_n$ of finite subsets such that $F_0=\{e\}$, $F_n=F_n^{-1}$. In view of \cite{b23}, it suffices to find a cellular subset $Y$ of $X$ with no asymptotically isolated balls.

Since $X$ is not cellular, there exists $F\in[G]^{<\w}$ such that

$(1)$ for every $n\in\NN$, there is $x\in X$ such that 
$$B_X^\square(x,F)\setminus B_X(x,F_n)\neq\varnothing.$$
We assume that $G$ is finitely generated and choose a system of generators $K\in[G]^{<\w}$ such that $K=K^{-1}$ and $F\subseteq K$. Then we consider the Cayley graph $\Gamma=Cay(G,K)$ with the set of vertices $G$ and the set of edges $\{\{g,h\}:g^{-1}h\in K\}$. We endow $\Gamma$ with the path metric $d$ and say that a sequence $a_0,...,a_n\in G$ is a geodesic path if $a_0,...,a_n$ is the shortest path from $a_0$ to $a_n$, in particular, $d(a_0,a_n)=n$. Using $(1)$, for each $n\in\NN$, we choose a geodesic path $L_n$ of length $3^n$ such that $L_n\subset X$ and 

$(2)$ $B_G(L_n,F_n\cap B_G(L_{n+1},F_{n+1})=\varnothing$ for every $n\in\NN$.\\ Let $L_n=\{a_{n0},...,a_{n3^n}\}$. For each $i\in\{0,...,3^n\}$, we take a tercimal decomposition of $i$ and denote by $Y_n$ the subset of all $a_i\in L$ such that $i$ has no $1$-s in its decomposition based on $\{0,1,2\}$ (see \cite{b26}). By \cite{b2} and the construction of $Y_n$, the set $Y=\bigcup_{n\in\NN}Y_n$ is cellular and has no asymptotically isolated balls.

Now let $G$ be an arbitrary countable group. We consider a subgroup $H$ of $G$ generated by $F$ and decompose $G$ into left cosets by $H$. If $X$ meets only finite number of these cosets then $X$ is contained in some finitely generated subgroup of $G$ and we arrive in the previous case. At last, let $\{Hx_n:n\in\NN\}$ be a decomposition of $G$ into left cosets by $H$ and $X$ meets infinitely many of them. We endow each $Hx_n$ with the structure of a graph $\Gamma_n$ naturally isomorphic to the Cayley graph $Cay(H,F)$. Then we use $(1)$ to choose an increasing sequence $(m_n)_{n\in\w}$ and a sequence $(L_n)_{n\in\NN}$ of geodesic paths of length $3^n$ satisfying $(2)$ and such that $L_n\subset X$. For each $n\in\NN$, we define a subset $Y_n$ of $L_n$ as before, and put $Y=\bigcup_{n\in\NN}Y_n$. By $(2)$ and the construction of $Y_n$, $Y$ is cellular an dhas no asymptotically isolated balls.

{\em Proof of Theorem~\ref{t1}}.

We follow the tour $(i)\Rightarrow(iv)\Rightarrow(ii)\Rightarrow(iii)\Rightarrow(i)$.

$(i)\Rightarrow(iv)$. We prove that a piecewise shifted $IP$-subset 
$$A=\{g_{i_1}g_{i_2}...g_{i_n}b_{i_n}:0\le i_1<...<i_n<\w\}$$
of $G$ is not scattered. For each $m\in\w$, let $$A_m=\{g_{i_1}g_{i_2}...g_{i_n}b_{i_n}:m< i_1<...<i_n<\w\}.$$
We take an arbitrary $p\in A^*$ and show that $\vt_p(A)$ is infinite.

If $A_m\in p$ for every $m\in\w$ then $g_np\in A^*$ for each $n\in\w$. Otherwise, there exists $m\in\w$ such that 
$$\{g_mg_{i_1}...g_{i_n}b_{i_n}:m< i_1<...<i_n<\w\}\in p.$$
Then $g_m^{-1}p\in A^*$ and we repeat the arguments for $g_m^{-1}p$.

$(iv)\Rightarrow(ii)$. Assume that $A$ is not disparse and take $p\in A^*$ such that $p$ is not isolated in $\vt_p(A)$. Then $p=qp$ for some $q\in G^*$. The set $\{x\in G^*:xp=p\}$ is a closed subsemigroup of $G^*$ and, by \cite[Theorem 2.5]{b1}, there is an idempotent $r\in G^*$ such that $p=rp$. We take $R\in r$ and $P_g\in p$, $g\in R$ such that $\bigcup_{g\in R}gP_g\subseteq A$.Since $r$ is an idempotent, by \cite[Theorem 5.8]{b1}, there is an injective sequence $(g_n)_{n\in\w}$ in $G$ such that $$\{g_{i_1}...g_{i_n}:0\le i_1<...<i_n<\w\}\subseteq R$$.
For each $n\in\w$, we pick $b_n\in\bigcap\{P_g: g=g_{i_1}...g_{i_n}:0\le i_1<...<i_n<\w\}$ and note that 
$$\{g_{i_1}...g_{i_n}b_{i_n}:0\le i_1<...<i_n<\w\}\subseteq A.$$

$(ii)\Rightarrow(iii)$. We assume that $A$ contains a subset coarsely equivalent to the group $B=\oplus_\w\ZZ_2$. Then there exist a subset $X$ of $B$, $H\in[B]^{<\w}$ such that $B=H+X$, and an injective $\prec$-mapping $f:X\to A$. We take an arbitrary idempotent $r\in B^*$, pick $h\in H$ such that $h+X\in r$ and put $p=r-h$. Since $r+p=r$, we see that $p$ is not isolated in $\vt_p(X)$. We denote $q=f^\beta(p)$. Let $b\in B$, $b\neq0$ and $b+p\in X^*$. Since $f$ is an injective $\prec$-mapping, there is $g\in G\setminus\{e\}$ such that $f^\beta(b+p)=g+q$. It follows that $q$ is not isolated in $\vt_q(A)$. Hence $A$ is not disparse.

$(iii)\Rightarrow(i)$. Let $X$ be a countable subset of $A$. By Proposition~\ref{p3}, $X$ is cellular. By \cite{b23}, $X$ satisfies Proposition~\ref{p1}$(ii)$. Hence $X$ is scattered. By Proposition~\ref{p2}, $A$ is scattered.\\\\
{\em Proof of Theorem~\ref{t2}}.

We assume that $\mu(A)>0$ for some Banach measure $\mu$ on $G$. We use the arguments from \cite[p. 506-507]{b10} to choose a decreasing sequence $(A_n)_{n\in\w}$ of subsets of $G$ and an injective sequence $(g_n)_{n\in\w}$ in $G$ such that $A_0=A$ and $\mu(A_n)>0$ $g_nA_{n+1}\subseteq A_n$ for each $n\in\w$. We pick $x_n\in A_{n+1}$ and put 
$$X=\{g_0^{\varepsilon_0}...g_n^{\varepsilon_n}x_n:n\in\w\text{ , } \varepsilon_i\in\{0,1\}\}.$$
By the construction $X$ is a piecewise shifted $IP$-sets and $X\subseteq A$. By Theorem~\ref{t1}, $X$ is not scattered.

\end{document}